\input amstex
\documentstyle{amsppt}
\pageheight{204mm}
\pagewidth{133mm}
\magnification\magstep1

%%%%% \amsusui.tex %%%%%

\def\nologo{\let\logo@\empty}

\def\Ad{\operatorname{Ad}}

\def\Aut{\operatorname{Aut}}
\def\Aut{\operatorname{Aut}}

\def\gr{\operatorname{gr}}
\def\Hom{\operatorname{Hom}}

\def\Im{\operatorname{Im}}

\def\Isom{\operatorname{Isom}}

\def\LMH{\operatorname{LMH}}

\def\Mor{\operatorname{Mor}}

\def\rank{\operatorname{rank}}

\def\SL{\operatorname{SL}}

\def\Spec{\operatorname{Spec}}

\def\toric{\operatorname{toric}}
\def\abtoric{\operatorname{|toric|}}
\def\torus{\operatorname{torus}}

\def\ts{\tsize\sum}

\def\bC{\bold C}

\def\be{\bold e}

\def\bN{\bold N}

\def\bQ{\bold Q}

\def\bR{\bold R}

\def\bZ{\bold Z}

\def\cA{{\Cal A}}
\def\cB{{\Cal B}}

\def\cO{{\Cal O}}

\def\cS{{\Cal S}}

\def\fg{{\frak g}}

\def\g{\gamma}
\def\G{\Gamma}

\def\sig{\sigma}
\def\Sig{\Sigma}

\def\vf{\varphi}

\def\.{$.\;$}
\def\an{{\text{\rm an}}}

\def\gp{{\text{\rm gp}}}
\def\loga{{\text{\rm log}}}
\def\mult{{\text{\rm mult}}}

\def\resp.{\text{\rm resp}.\;}

\def\O^logten{\cO\^log\otimes}

\let\bs=\backslash

\let\lan=\langle
\let\lan=\langle

\let\hra=\hookrightarrow

\let\ran=\rangle

\let\sub=\subset

\let\x=\times
\let\v=\vee

\def\Dc{\check{D}}
\def\Ec{\check{E}}

%%%%%

\topmatter

\title
Moduli of log
mixed Hodge structures 
\endtitle

\author
Kazuya Kato, Chikara Nakayama, Sampei Usui
\endauthor

\address
\newline
{\rm Kazuya KATO}
\newline
Department of Mathematics
\newline
University of Chicago
\newline
5734 S.\ University Avenue
\newline
Chicago, Illinois 60637, USA
\newline
{\tt kkato\@math.uchicago.edu}
\endaddress

\address
\newline
{\rm Chikara NAKAYAMA}
\newline
Graduate School of Science and Engineering
\newline
Tokyo Institute of Technology 
\newline
Meguro-ku, Tokyo, 152-8551, Japan
\newline
{\tt cnakayam\@math.titech.ac.jp}
\endaddress

\address
\newline
{\rm Sampei USUI}
\newline
Graduate School of Science
\newline
Osaka University
\newline
Toyonaka, Osaka, 560-0043, Japan
\newline
{\tt usui\@math.sci.osaka-u.ac.jp}
\endaddress

\abstract
We announce the construction of toroidal partial compactifications 
of the moduli spaces of mixed Hodge structures with polarized graded quotients.
They are moduli spaces of log mixed Hodge structures with polarized graded 
quotients.  
We include an 
application to the analyticity of zero loci of some functions. 
\endabstract

\endtopmatter

\document

\head
\S0. Introduction
\endhead

Log Hodge structure is a natural formulation of \lq\lq degenerating 
family of Hodge
structures".
In \cite{KU09}, the moduli spaces of polarized log Hodge
structures were constructed. 
In this paper, we construct the moduli spaces of log mixed Hodge structures whose graded quotients by weight filtrations are polarized log Hodge
structures.
The construction is parallel to the pure case \cite{KU09}. 
We add points at infinity to the non-log moduli by using the mixed version of 
nilpotent orbits.  
It is also parallel to the pure case to prove that the constructed spaces are
actually the moduli of log mixed Hodges structures with polarized
graded quotients.
As in the pure case, they are like toroidal partial compactifications with slits of the moduli of mixed Hodge structures with polarized graded quotients. 

We omit here the details of proofs of the above facts, which are to be 
published in the series of papers \cite{KNU.p1}, \cite{KNU.p2}, ....

In the final section, we include an application on the analyticity of zero loci of some functions to show how our spaces are helpful in studying such geometric problems.

\vskip20pt

\head
\S1. Mixed Hodge structures and moduli
\endhead

We review the construction of the moduli spaces of mixed Hodge structures with polarized graded quotients, studied in \cite{U84}. 
This is a mixed Hodge theoretic version of Griffiths domains in \cite{G68}.

\medskip

{\bf 1.1.} 
We fix a $4$-ple
$$
\Lambda=(H_0, W, (\lan\;,\;\ran_k)_k, ((h^{p,q}_k)_{p,q})_k),
$$
where 

$H_0$ is a finitely generated free $\bZ$-module, 

$W$ is an increasing filtration on $H_{0,\bQ}:=\bQ\otimes_\bZ H_0$, 

$\lan\;,\;\ran_k$ is a non-degenerate $\bQ$-bilinear form $\gr^W_k \times
\gr^W_k \to \bQ$ given for each $k\in \bZ$ which is symmetric if $k$ is even and anti-symmetric if $k$ is odd, and 

$h^{p,q}_k$ is a non-negative integer given for $p, q, k\in\bZ$ such that $h^{p,q}_k = 0$ unless $p + q = k$, 
that $h^{p,q}_k = h^{q,p}_k$ for all $p, q, k$, and that 
$$
\rank_\bZ(H_0) = \tsize\sum_{p,q,k} h^{p,q}_k, \quad
\dim_\bQ(\gr^W_k)= \tsize\sum_{p,q} h^{p,q}_k \;\;\text{for all} \;k.
$$

\medskip

{\bf 1.2.} 
We fix notation.

For $A=\bZ, \bQ, \bR,$ or $\bC$, let
$G_A$ be the group of all $A$-automorphisms of $H_{0,A}$ which preserve
$A \otimes_\bZ W$ and $A \otimes_\bZ \lan\;,\ran_k$ for any $k$. 

For $A=\bQ, \bR, \bC$, let
$\fg_A$ be the set of all $A$-homomorphisms $N : H_{0,A}\to H_{0,A}$
satisfying the following conditions (1) and (2). 

\medskip 

(1) $N(A \otimes_\bQ W_k) \subset A \otimes_\bQ W_k$ for any $k$.

\medskip

(2) For any $k$, the homomorphism
$\gr^W_k(N) : A \otimes_\bQ \gr^W_k \to A \otimes_\bQ \gr^W_k$ induced by $N$ satisfies $\lan \gr^W_k(N)(x),
y\ran_k+ \lan x, \gr^W_k(N)(y)\ran_k=0$ for all $x, y\in A \otimes_\bQ \gr^W_k$.
\medskip

{\bf 1.3.}
Let $D$ be the set of all decreasing filtration $F$ on $H_{0,\bC}$ for which $(H_0, W, F)$ is a mixed Hodge structure such that the $(p, q)$ Hodge number of $F(\gr^W_k)$
coincides with $h^{p,q}_k$ for any $p, q, k\in \bZ$ and such that $F(\gr^W_{k})$ is polarized by $\lan\;,\;\ran_k$ for all $k$.

Let $\Dc \supset D$ be the set of all decreasing filtrations $F$ on $H_{0,\bC}$ 
satisfying the following properties (1) and (2).

\medskip

(1) $\dim \gr^p_F (\gr^W_{k,\bC}) = h^{p,q}_k$ for any $p, q, k\in \bZ$ such that $p+q=k$.

\medskip

(2) $\lan\;,\;\ran_k$ kills $F^p(\gr^W_{k})\times F^q(\gr^W_{k})$ for any $p,q,k\in \bZ$ such that $p+q>k$.

\medskip

Then $G_{\bC}= \Aut(H_{0,\bC}, W, (\lan\;,\;\ran_k)_k)$ acts transitively on 
$\Dc$ and hence $\Dc$ is an analytic manifold. 
Furthermore, $D$ is open in $\Dc$ and hence it is also an analytic manifold.

\medskip
This $D$ was introduced in \cite{U84} and is 
a natural generalization of Griffiths classifying space of polarized
Hodge structures (\cite{G68}) 
which is the case where there is only one $k$ such that $\gr^W_k\neq 0$.

\vskip20pt

\head
\S2. Space of nilpotent orbits
\endhead

\head
\S2.1.  Set $D_\Sig$.
\endhead

 We fix $\Lambda=(H_0, W, (\lan\;,\;\ran_k)_k, (h^{p,q}_k)_k)$ as in 1.1.

\medskip

{\bf 2.1.1.} Nilpotent cone.

A finitely generated and sharp cone of $\fg_{\bR}$ is called a {\it nilpotent cone} if 
it is generated by mutually commuting nilpotent elements.

Note that, since a nilpotent cone is finitely generated, there are only finitely many faces of it.

\medskip

{\bf 2.1.2.} 
Fan in $\fg_\bQ$.

A {\it fan $\Sig$ in $\fg_{\bQ}$} is a set of nilpotent cones in $\fg_\bR$ satisfying the following conditions (1)--(4).

\medskip

(1) All $\sig \in \Sig$ are rational. 
That is, $\sig$ is  generated by (a finite number of) elements of $\fg_\bQ$ over $\bR_{\ge0}$. 

\medskip

(2) If $\sig\in \Sig$, all faces of $\sig$ belong to $\Sig$.

\medskip

(3) If $\sig, \sig' \in \Sig$, $\sig \cap \sig'$ is a face of $\sig$.

\medskip

(4) (Admissibility) 
Let $\sig\in \Sig$. 
Then for any element $N$ of the interior of $\sig$, the relative monodromy filtration $M(N, W)$ exists. 
Furthermore, this filtration is independent of the choice of $N$.

\medskip

{\bf 2.1.3.} 
Nilpotent orbit. 

Let $D$ and $\Dc$ be the spaces in 1.3.

Let $\sig$ be a nilpotent cone.
A subset $Z$ of $\Dc$ is called a {\it $\sigma$-nilpotent orbit} (resp\. {\it $\sig$-nilpotent $i$-orbit}) if the following conditions (1)--(3) are satisfied for some $F\in Z$.

\medskip

(1) $Z=\exp(\sig_\bC)F$ (resp\. $Z=\exp(i\sig_\bR)F$). Here $\sig_\bC$ (resp\. $\sig_\bR$) is the
$\bC$ (resp\. $\bR$)-linear span of $\sig$ in $\fg_\bC$ (resp\. $\fg_\bR$).

\medskip

(2) $N(F^p)\subset F^{p-1}$ for all $p\in \bZ$.

\medskip

(3) Take a finite family $(N_j)_{1\leq j\leq n}$ of elements of $\sig$ which 
generates $\sig$. 
Then, if $y_j \in \bR$ and $y_j$ are sufficiently large for $1\le j \le n$, we have $\exp(\sum_{j=1}^n iy_jN_j)F\in D$.

\medskip

Note that, if (1)--(3) are satisfied for some $F\in Z$, they are satisfied for all $F\in Z$.

\medskip

{\bf 2.1.4.} 
Let $D_\Sig$ (resp\.$D^{\sharp}_\Sig$) be the set of all pairs $(\sig, Z)$, where
$\sig\in \Sig$ and $Z$ is a $\sig$-nilpotent orbit (resp\. $\sig$-nilpotent $i$-orbit). 

\medskip

We have embeddings 
$$
D \subset D_\Sig \quad \text{and} \quad
D \subset D_\Sig^{\sharp}, \quad 
F \mapsto (\{0\}, \{F\}).
$$

We have a surjection
$$
D_\Sig^\sharp \to D_\Sig, \quad
(\sig, Z) \mapsto (\sig, \exp(\sig_\bC)Z).
$$

{\bf 2.1.5.} Compatibility with $\G$.
Let $\G$ be a subgroup of $G_{\bZ}$. 
We say $\Sig$ and $\G$ are {\it compatible}, if for any $\g \in \G$ and $\sig \in \Sig$, we have $\Ad(\g)\sig \in \Sig$. 

Then $\Gamma$ acts on $D_\Sig$ and also on $D_{\Sig}^{\sharp}$ by $(\sig, Z)\mapsto (\Ad(\gamma)\sig, \gamma Z)$ ($\gamma\in \Gamma$). 

Further, we say $\Sig$ and $\G$ are {\it strongly compatible} if they are compatible and for any $\sig \in \Sig$, any element of $\sig$ can be written as a finite sum of elements of the form $aN$, where $a \in \bR_{\ge0}$ and $N \in \sig$ satisfying $\exp(N) \in \G$. 

\medskip

\head
\S2.2. Set $E_\sig$
\endhead

  Let $\Sig$ and $\Gamma$ be as in \S2.1. 
  Assume that they are strongly compatible. 
  We fix $\sig \in \Sig$ in this \S2.2. 

\medskip

{\bf 2.2.1.}
Let $D_\sig=D_{\text{face}(\sig)}$, $D_\sig^\sharp=D_{\text{face}(\sig)}^\sharp$, where $\text{face}(\sig)$ denotes the fan consisting of all faces of $\sig$.

Let $\G(\sig)=\G \cap \exp(\sig) $.   
It is a sharp and torsion free fs monoid (for the terminology, see \cite{K89}, \cite{KU09}). 
The associated group $\G(\sig)^\gp$ is a finitely generated abelian group and it is strongly compatible with the fan $\text{face}(\sig)$.

We will 
regard $\G(\sig)^\gp\bs D_\sig$ as a quotient of a subset $E_\sig$ of an analytic space ${\check E}_{\sig}$ explained below 
(see \cite{KU09} for the pure version). 

\medskip

{\bf 2.2.2.}
Associated to $\sig$, we have the torus and the toric variety: 
$$
\torus_\sig:=\Spec(\bC[\G(\sig)^{\v\,\gp}])_\an \sub 
\toric_\sig:=\Spec(\bC[\G(\sig)^\v])_\an.
$$
  Here and hereafter we denote $\Hom(P,\bN)$ by $P^\v$ for an 
fs monoid $P$. 

We denote
$$
\Ec_\sig = \toric_\sig\x \Dc.
$$

We use some facts about log structures in \cite{K89}, \cite{KU09}.
We endow $\toric_\sig$ with the canonical fs log structure, and endow $\Ec_\sig$ with its inverse image.
Since $\Dc$ is smooth, $\Ec_\sig$ is a logarithmically smooth fs log analytic space.

Let $(\toric_\sig^\loga, 
\cO_{\toric_\sig}^\loga)$ be the associated ringed space.
We have isomorphisms 
$$
\pi_1(\toric_\sig^\loga) \simeq \pi_1(\torus_\sig)\simeq \G(\sig)^\gp.
$$

\medskip

{\bf 2.2.3.}
Let $q \in \toric_\sig$.
Let $\sig(q)$ be the face of $\sig$ corresponding to the face $\pi_1^+(q^\loga):= \pi_1(q^\loga) \cap \Gamma(\sig)$ of $\G(\sig)$.
Let $\cS'=\{f\in \G(\sig)^\v\;|\;f(q)\neq 0\}$, where $f$ is regarded as a holomorphic function on $\toric_\sig$.
Let $\be: \sig_\bC/(\sig(q)_\bC+\log(\G(\sig)^\gp)) @>\sim>>\;\Hom((\cS')^\gp, \bC^\x)$ be the isomorphism defined by $(\be(z\log\g), f) = \exp(2\pi iz(\g, f))$ for $z\in \bC$, $\g \in \G(\sig)^\gp$, $f \in (\cS')^\gp$.
Let $z$ be an element of $\sig_\bC$ whose image in $\sig_\bC/(\sig(q)_\bC+\log(\G(\sig)^\gp))$ coincides with the class of $q\in \Hom((\cS')^\gp,  \bC^\x)$ under the above isomorphism.

We define the subset $E_\sig$ of $\Ec_\sig$ by the following condition.

\medskip

For $(q, F) \in \Ec_\sig= \toric_\sig \times \Dc$, $(q, F)\in E_\sig$ if and only if $\exp(\sig(q)_\bC)\exp(z)F$ is a $\sig(q)$-nilpotent orbit.

\medskip

Denote
$$
\align
&\abtoric_\sig
=\Hom(\G(\sig)^\v,\bR_{\ge0}^\mult)\sub \toric_\sig=\Hom(\G(\sig)^\v,\bC^\mult),\endalign
$$
where $\bR_{\ge0}^\mult$ and $\bC^\mult$ are the sets $\bR_{\ge0}$ and $\bC$ regarded as monoids by multiplication, respectively.
Define
$$
\Ec_\sig^\sharp=\abtoric_\sig \x \Dc \sub \Ec_\sig,\quad
E_\sig^\sharp=E_\sig \cap \Ec_\sig^\sharp.
$$

Then the subset $E_\sig^\sharp$ of $\Ec_\sig^\sharp$ can be characterized by the following condition.

\medskip

For $(q, F) \in \Ec_\sig^\sharp= \abtoric_\sig \times \Dc$, $(q, F)\in E_\sig^\sharp$ if and only if $\exp(i\sig(q)_\bR)\exp(iy)F$ is a $\sig(q)$-nilpotent $i$-orbit.
Here $y\in \bR$ is the imaginary part of the above $z$.

\medskip

{\bf 2.2.4.}
Define canonical maps $\vf : E_\sig \to \G(\sig)^\gp\bs D_\sig$ and $\vf^\sharp : E^\sharp_\sig \to D^\sharp_\sig$ by
$$
\vf(q, F)=((\sig(q), \exp(\sig(q)_\bC)\exp(z)F) \bmod \G(\sig)^\gp),
$$
$$
\vf^\sharp(q,F)=(\sig(q),\exp(i\sig(q)_\bR)\exp(iy)F).
$$

\medskip

\head 
\S2.3. Strong topology, the category $\cB(\log)$, log manifolds
\endhead

In \S2.4 below, we will endow $\G\bs D_{\Sig}$  with a structure of a local ringed space over $\bC$ with an fs log structure, and we define a topology on $D^\sharp_\Sig$. 
In this \S2.3, we give preparation for them.

\medskip

{\bf 2.3.1.}  
Strong topology.

Let $Z$ be an analytic space, and $S$ a subset. 
The {\it strong topology} of $S$ in $Z$ is defined as follows. A subset $U$ of $S$ is open for this topology if and only if for any analytic space $A$ and any morphism $f : A \to Z$ of analytic spaces such that $f(A)\subset S$, $f^{-1}(U)$ is open in $A$. 
It is stronger than or equal to the topology as a subspace of $Z$.

\medskip

{\bf 2.3.2.} 
The category $\cB(\log)$.

As in \cite{KU09}, in the theory of moduli spaces of log mixed 
Hodge structures, we have to enlarge the category of analytic spaces because the moduli spaces are often not analytic spaces.

Let $\cA$ be the category of anaytic spaces and let $\cA(\log)$ be the category of fs log analytic spaces. 
We enlarge $\cA$ and $\cA(\log)$ to $\cB$ and $\cB(\log)$, respectively, as
follows. In \cite{KU09}, the moduli spaces live in $\cB(\log)$. 

$\cB$ (resp\. $\cB(\log)$) is the category of all local ringed spaces $S$ over
$\bC$ (resp\. local ringed spaces $S$ 
over $\bC$ endowed with an fs log structure) having the following property:
$S$ is locally isomorphic to a subset of an analytic space (resp\. an fs log analytic space) $Z$ with the strong topology in $Z$ (2.3.1) and with the inverse image (resp.\ inverse images) of the sheaf of rings $\cO_Z$ (resp\. $\cO_Z$ and the log structure $M_Z$).

\medskip

{\bf 2.3.3.} 
A log manifold.

A {\it log manifold} is an object $S$ of $\cB(\log)$ such that locally on $S$, 
we can take a logarithmically smooth 
fs log analytic space $Z$ 
(for the terminology, see, e.g., \cite{KU09}) and $\omega_1,\cdots, \omega_n\in \G(Z, \omega^1_{Z})$ ($\omega^1_{Z}$ is the sheaf of differential forms with 
log poles; see \cite{KU09}) such that $S$ (as an object of $\cB(\log)$) is isomorphic to an open subspace of the subspace $\{z\in Z\;|\;\text{$\omega_1, \cdots, \omega_n$ are zero in $\omega^1_z$}\}$ with the strong topology in $Z$ 
and the inverse images of $\cO_Z$ and $M_Z$. 

In \cite{KU09}, we saw that moduli spaces of polarized log Hodge structures were log manifolds.

\medskip

\head
\S2.4. Topology, local ringed space structure, log structure of $\G\bs D_\Sig$
\endhead

Let $\Sigma$ be a fan in $\fg_{\bQ}$ (2.1.2). 
Let $\Gamma$ be a subgroup of $G_\bZ$, which is strongly compatible with $\Sig$. 
Let $\sig \in \Sig$. 

\medskip

{\bf 2.4.1.}
We endow the subset $E_\sig$ of $\Ec_\sig$ in \S 2.2 with the following structures of log local ringed spaces over $\bC$.
The topology is the strong topology in $\Ec_\sig$.
The sheaf $\cO$ of rings and the log structure $M$ are the inverse images of $\cO$ and $M$ of $\Ec_\sig$, respectively.

We endow $\G\bs D_\Sig$ with the strongest topology for which
the maps 
$\pi_\sig:E_\sig \overset \vf \to \to\G(\sig)^\gp\bs D_\sig\to\G\bs D_\Sig$
are continuous
for all $\sig\in\Sig$.
  Here $\vf$ is as in 2.2.4.
We endow $\G\bs D_\Sig$ with the following sheaf
of rings $\cO_{\G\bs D_\Sig}$ over $\bC$ and the
following log structure $M_{\G\bs D_\Sig}$.
For any open set $U$ of $\G\bs D_\Sig$ and for any
$\sig\in\Sig$, let $U_\sig:=\pi_\sig^{-1}(U)$ and
define
$$
\align
&\cO_{\G\bs D_\Sig}(U)\;
\text{(\resp. $M_{\G\bs D_\Sig}(U)$)}\\
&:=\{\text{map $f:U\to\bC$}\;|\;
f\circ\pi_\sig\in\cO_{E_\sig}(U_\sig)\;
(\resp. \in M_{E_\sig}(U_\sig))\;
(\forall\sig\in\Sig)\}.
\endalign
$$

\medskip

{\bf 2.4.2.}
We introduce the topology of $E_\sig^\sharp$ as a
subspace of $E_\sig$ (2.2.3).
We introduce on $D_\Sig^\sharp$ the strongest
topology for which the maps
$E_\sig^\sharp \overset {\vf^\sharp} \to \to 
D_\sig^\sharp\to D_\Sig^\sharp$ $(\sig\in\Sig)$
are continuous.
Here $\vf^\sharp$ is as in 2.2.4. 
Note that the surjection $D_\Sig^\sharp
\to\G\bs D_\Sig$ (cf.\ 2.1.4) becomes continuous.

\medskip

{\bf 2.4.3.}
The above topologies have the following properties.

\medskip

(1) Let $(\sig,Z)\in D_\Sig$, let $F\in Z$, and
write $\sig=\ts_{1\le j\le n}\bR_{\ge0}N_j$.
Then
$$
\big((\sig,Z)\bmod\G\big)
=\lim\Sb\Im(z_j)\to\infty\\1\le j\le n\endSb
\big(\exp\big(\ts_{1\le j\le n}z_jN_j\big)F
\bmod\G\big)\quad\text{in $\G\bs D_\Sig$}.
$$

(2) Let $(\sig,Z)\in D_\Sig^\sharp$, let
$F\in Z$, and let $N_j$ be as above.
Then
$$
(\sig,Z)=\lim\Sb y_j\to\infty\\1\le j\le n\endSb
\exp\big(\ts_{1\le j\le n}iy_jN_j\big)F\quad
\text{in $D_\Sig^\sharp$}.
$$

\medskip

\head
\S2.5. Main theorem A
\endhead

\medskip

\proclaim{Theorem A}  
Let $\Sigma$ be a fan in $\fg_{\bQ}$ $(2.1.2)$. 
Let $\Gamma$ be a subgroup of $G_{\bZ}$, which is strongly compatible with $\Sig$. 
Let $\sig \in \Sig$. 
\smallskip

$(1)$ $E_\sig$ belongs to $\cB(\log)$. It is a log manifold.
\smallskip

$(2)$ $E_\sig\to \G(\sig)^\gp \bs D_\sig$ is a $\sig_\bC$-torsor in the category $\cB(\log)$. 
$E^\sharp_\sig \to D^\sharp_\sig$ is an $i\sig_\bR$-torsor in the category of topological spaces. 
\smallskip

$(3)$  The action of $\Gamma$ on $D^{\sharp}_\Sig$ is proper. 
The spaces $\G \bs D^{\sharp}_\Sig$ 
and $\G \bs D_{\Sig}$ are Hausdorff.
\smallskip

$(4)$ Assume that $\G$ is neat.  
Then $D^\sharp_{\Sig} \to \Gamma \bs D_\Sig^{\sharp}$ is a local homeomorphism.
\smallskip

$(5)$ Assume that $\G$ is neat.  
Then $\Gamma \bs D_\Sig$ belongs to $\cB(\log)$. 
It is a log manifold and $(\G\bs D_\Sig)^\loga= \G \bs D_\Sig^\sharp$. 
\endproclaim

\vskip20pt

\head
\S3. Log mixed Hodge structures and moduli
\endhead

Let $S$ be an object in $\cB(\log)$.
Then we have the associated ringed space $(S^\loga, {\cO}_S^\loga)$. 
See \cite{KU09} \S2.2 for the definition.

\head
\S3.1. Log mixed Hodge structures
\endhead

{\bf 3.1.1.} 
Polarized log Hodge structure.

A polarized log Hodge structure was the main ingredient of \cite{KU09}. 
See \S2.4 in loc.\ cit.\ for the definition. 

The definition includes two points. 

(1) PHS after sufficiently twisted specializations.

(2) Small Griffiths transversality.

\medskip

{\bf 3.1.2.} 
For example, degeneration of elliptic curves gives a polarizable pure log Hodge structure of weight $1$ (not mixed!).
This is explained in detail in \cite{KU09} \S0. 

\medskip

{\bf 3.1.3.} 
Log mixed Hodge structure with polarized graded quotients.

This is the main ingredient of this paper. 
Let $S$ be an object in $\cB(\log)$. 
A {\it log mixed Hodge structure with polarized graded quotients} ({\it LMH with PGQ,} for short) over $S$ is a $4$-ple $(H_\bZ, W, (\langle\;,\;\rangle_k)_k, 
F)$, where

$H_\bZ$ is a locally constant sheaf of finitely generated free $\bZ$-modules on $S^\loga$,

$W$ is an increasing filtration on $H_\bQ:=\bQ \otimes_\bZ H_\bZ$,

$\langle\;,\;\rangle_k$ is a $(-1)^k$-symmetric bilinear form $\gr^W_k \times \gr^W_k \to \bQ$ given for each $k\in \bZ$, 

$F$ is a decreasing filtration of the $\cO_S^\loga$-module 
$\cO_S^\loga \otimes_\bZ H_\bZ$, 

\noindent 
such that $(H_{\bZ}, W_{\bR}, F)$ is an LMH in the sense of \cite{KU09} \S2.6 
and such that for each $k \in \bZ$, the induced data on $\gr^W_k$ form a polarized log Hodge structure of pure weight $k$.

  Thus there are three points in the definition. 

(1) MHS with polarized graded quotients 
after sufficiently twisted specializations. 

(2) Small Griffiths transversality.

(3) Admissibility of local monodromy.

\medskip

{\bf 3.1.4.}
The key observation in the pure case in \cite{KU09} 0.4.25 can be generalized, and we have also in the present case the following:

\medskip

(an LMH with PGQ on an fs log point)=
(a nilpotent orbit in the mixed case).

\medskip

\head
\S3.2. Moduli functor 
\endhead

We define the moduli functor of log mixed Hodge structure with polarized graded
quotients.

\medskip

{\bf 3.2.1.} 
Fix $\Phi= (\Lambda, \Sig, \G)$, where $\Lambda$ is as in 1.1, $\Sig$ and $\G$ are as in \S2.1, and $\Sig$ is assumed to be strongly compatible with $\G$. 

\medskip

{\bf 3.2.2.} 
Let $S$ be an object of $\cB(\log)$. 
By {\it a log mixed Hodge structure of type $\Phi$} over $S$, we mean an LMH 
with polarized graded quotients 
$H=(H_\bZ, W, (\langle\;,\;\rangle_k)_k), 
F)$ endowed with a global section $\mu$ of the sheaf $\Gamma \bs\Isom((H_\bZ, W, (\langle\,,\,\rangle_k)_k), (H_0, W, (\langle\,,\,\rangle_k)_k))$ on $S^\loga$ which satisfies the following conditions (1) and (2). 

\medskip

(1) $\rank_\bZ(H_\bZ)=\sum_{p,q,k} h^{p,q}_k, \quad \rank_{\cO_S^{\log}} 
(F^p)= \sum_{k\in\bZ, r\geq p} h_k^{r, k-r}\quad \text{for all} \, p$.

\medskip

(2) For any $s\in S$ and $t \in S^\loga$ lying over $s$, if ${\tilde \mu}_t :
(H_{\bZ,t}, W, (\langle\;,\;\rangle_k)_k)\overset{\simeq}\to \to
(H_0, W,(\langle\;,\;\rangle_k)_k)$ is a representative of the stalk of $\mu$ at $t$, then
there exists $\sig\in \Sig$ such that the image of the composite map
$$
\align
\Hom(M_{S,s}/\cO_{S,s}^\times, \bN) &\hra \pi_1(\tau^{-1}(s)) \\
&\to\Aut(H_{\bZ,t}, W, (\langle\;,\;\rangle_k)_k)@>{\text{by} \;{\tilde \mu}_t}>> \Aut(H_0,
W, (\langle\;,\;\rangle_k)_k)
\endalign
$$
is contained in $\exp(\sig)$. 
  Furthermore, if we take the smallest such $\sig \in \Sig$,  
then the $\exp(\sig_{\bC})$-orbit $Z$ including 
$\tilde \mu_t(\bC\otimes_{{\cO}^{\log}_{S,t}}F_t )$, which is independent of the choice of a $\bC$-algebra homomorphism ${\cO}^{\log}_{S,t} \to \bC$, 
is a $\sig$-nilpotent orbit (cf. 3.1.4, cf. also [KU09] 0.4.24, 2.5.1, 2.5.5).

\medskip

We call a log mixed Hodge structure of type $\Phi$ over $S$ also a {\it log mixed Hodge structure with polarized graded quotients, with global monodromy in $\Gamma$, and with local monodromy in $\Sig$}. 

\medskip

{\bf 3.2.3.}  
Let $\LMH_{\Phi} : \cB(\log) \to$ (set) be the contravariant functor defined 
as follows: 
$\LMH_{\Phi}(S)$ for an object $S$ of $\cB(\log)$ is the set of isomorphism classes of log mixed Hodge structures of type $\Phi$ over $S$. 

\medskip

\head
\S3.3. Main theorem B
\endhead

\medskip

\proclaim{Theorem B} 
Assume that $\G$ is neat. 
Then the functor $\LMH_\Phi$ in $3.2.3$ is represented by $\G \bs D_\Sig$. 
\endproclaim

The period map $\LMH_\Phi\to \Mor(\bullet, \Gamma \bs D_\Sig)$ which is the isomorphism in this Theorem B is as follows. 
Let $S$ and $H$ be as in 3.2.2. 
Let $s \in S$. 
The associated point of $\Gamma \bs D_\Sig$ by this period map is the image 
of $(\sig,Z)\in D_{\Sig}$ in $\G \bs D_\Sig$, 
where $\sig$ and $Z$ are the ones in 
the last sentence of 3.2.2 (2). 

\medskip

The proofs of the theorems so far in this paper are similar to those in the pure case (\cite{KU09}). 
In the pure case, the key tool is the $\SL(2)$-orbit theorem in several variables of Cattani-Kaplan-Schmid \cite{CKS86}.
Instead of this, here we use a mixed Hodge theoretic version \cite{KNU08} 
of the $\SL(2)$-orbit theorem in several variables.

\vskip20pt

\head
\S4. Some application
\endhead

\proclaim{Theorem} 
Let $S$ be a complex analytic manifold, let $T$ be a smooth divisor on $S$, and let $S^*=S-T$. Let $\G$ be a neat subgroup of $G_\bZ$. 
For $1\leq j\leq n$, let $f_j: S^*\to \G \bs D$ be the 
period maps associated to variations of mixed Hodge structure $H_j$ 
with polarized graded quotients which are admissible with respect to $S$. 

Let $V=\{s\in S^*\;|\; f_1(s)=\dots=f_n(s)\}$, and let $\bar V$ be the closure of $V$ in $S$. 
Then $\bar V$ is an analytic subset of $S$.  
\endproclaim

{\it Proof.}
By a standard argument, we may assume that the
local monodromy of $H_j$ along the divisor $T$ is
unipotent. 
Then, by the admissibility of $f_j$, $H_j$ extends
to a log mixed Hodge structure $\tilde H_j$ over $S$ (see [KNU08] \S12).
Hence the period map $f_j$ of $H_j$ extends to a morphism 
$\bar f_j: S\to \G \bs D_{\Xi}$ in $\cB(\log)$ 
corresponding to $\tilde H_j$, 
where $\Xi$ is the fan consisting of all the one-dimensional rational nilpotent cones and $\{0\}$ (see [KU09] 4.3.1 (i) for the pure case). 
To prove the Theorem, 
it is sufficient to see that $\{s\in S\;|\;\bar f_1(s)=\dots =\bar f_n(s)\}$ is a closed analytic subset of $S$. Hence we are reduced to:

\proclaim{Proposition} 
Let $S$ be a complex analytic space, let $Y$ be an object of $\cB$, let $f_j: S\to Y$ ($1\leq j\leq n$) be morphisms in $\cB$, and assume that $Y$ is Hausdorff as a topological space. 
Let $C=\{s\in S\;|\; f_1(s)=\dots=f_n(s)\}$. 
Then $C$ is a closed analytic subset of $S$.
\endproclaim

{\it Proof of Proposition.}
Let $f=(f_j)_j: S\to Y^n$. 
Working locally on $Y$, we may assume that $Y$ is a subset of a complex analytic space $Z$ and $\cO_Y$ is the inverse image of $\cO_Z$.  
Then $C$ is the fiber product of $S\to Z^n\leftarrow Z$, where $S\to Z^n$ is the composite $S\to Y^n\to Z^n$ and $Z\to Z^n$ is the diagonal morphism. 
This proves Proposition.\qed

\medskip

{\it Remark 1.} 
Saito (\cite{S.p}) proved the case where $V$ is the zero locus of 
an admissible normal function. 
Note that our theorem is applied to the case of intermediate Jacobian (cf.\ \cite{KNU.s}).  
In particular, if we take $n=2$, an admissible normal function as $f_1$, 
and the zero section as $f_2$ in this case, then 
the set $V$ is nothing but the zero locus of $f_1$.

\medskip

{\it Remark 2.} 
In the series of papers \cite{KNU.p1}, \cite{KNU.p2}, ..., we plan to 
give a generalization of the Theorem in which $T$ can be any closed analytic 
subspace of $S$. 
This will give an alternative proof of a result of Brosnan and Pearlstein 
(\cite{BP.p}) when $V$ is the zero locus of an admissible normal function.


\Refs

\widestnumber\key{KNU.p2}

\ref
\key BP.p
\by P.\ Brosnan and G.\ Pearlstein
\paper On the algebraicity of the zero locus of an admissible normal function
\jour preprint
\vol 
\yr 
\pages 
\endref

\ref
\key CKS86
\by E.\ Cattani, A.\ Kaplan and W.\ Schmid
\paper Degeneration of Hodge structures
\jour Ann. of Math.
\vol 123
\yr 1986
\pages 457--535
\endref

\ref
\key G68
\by P.\ Griffiths
\paper Periods of integrals on algebraic manifolds. I.
Construction and properties of modular varieties
\jour Amer. J. Math.
\vol 90
\yr 1968
\pages 568--626
\endref

\ref
\key K89
\by K. Kato
\paper Logarithmic structures of Fontaine-Illusie
\inbook in \lq\lq Algebraic analysis, geometry, and number theory"
\ed J.-I. Igusa
\publ Johns Hopkins University Press
\publaddr Baltimore
\yr 1989
\pages 191--224
\endref

\ref
\key KNU08
\by K.\ Kato, C.\ Nakayama and S.\ Usui
\paper $\SL(2)$-orbit theorem for degeneration of mixed Hodge structure
\jour J.\ Algebraic Geometry
\vol 17
\yr 2008
\pages 401--479
\endref

\ref
\key KNU.s
\bysame
\paper Log intermediate Jacobians
\jour preprint, submitted
\vol 
\yr 
\pages 
\endref

\ref
\key KNU.p1
\bysame
\paper Classifying spaces of degenerating mixed Hodge structures, II\rom:
Spaces of $\SL(2)$-orbits
\jour in preparation
\yr
\pages
\endref

\ref
\key KNU.p2
\bysame
\paper Classifying spaces of degenerating mixed Hodge structures, III\rom:
Spaces of nilpotent orbits
\jour in preparation
\yr 
\pages 
\endref

\ref
\key KU09
\by K.\ Kato and S.\ Usui
\book Classifying spaces of degenerating polarized 
Hodge structures
\bookinfo  Ann.\ of Math.\ Stud., 
{\bf 169}
\publ Princeton Univ.\ Press
\publaddr  Princeton, NJ
\yr 2009
\endref

\ref
\key S.p
\by M.\ Saito
\paper  Hausdorff property of the N\'eron models of Green, Griffiths and Kerr
\jour arXiv:0803.2771
\vol 
\yr 
\pages 
\endref

\ref
\key U84
\by S.\ Usui
\paper Variation of mixed Hodge structure arising from
family of logarithmic deformations II\rom: Classifying space
\jour Duke Math\.J.
\vol 51-4
\yr 1984
\pages 851--875
\endref

\endRefs
\enddocument